\documentclass[a4paper, reqno]{amsart}
\usepackage{a4}
\usepackage{enumerate}
\usepackage[russian,english]{babel}
\usepackage{amsmath, amssymb,color,ulem}

\newtheorem{theorem}{Theorem}

\newtheorem{lemma}[theorem]{Lemma}

\newcommand{\ud}{\mathrm{d}}

\begin{document}

\title[Steady states in a structured epidemic model]{Steady states in a structured epidemic model with Wentzell boundary condition}

\author[\`{A}. Calsina]{\`{A}ngel Calsina}
\address{\`{A}ngel Calsina, Department of Mathematics, Universitat Aut\`{o}noma de Barcelona, Bellaterra, Spain}
\email{acalsina@mat.uab.es}

\author[J.\,Z. Farkas]{J\'{o}zsef Z. Farkas}
\address{J\'{o}zsef Z. Farkas, Institute of Computing Science and Mathematics, University of Stirling, Stirling, FK9 4LA, United Kingdom }
\email{jzf@maths.stir.ac.uk}

\subjclass{92D25, 47N60, 47D06, 35B35}
\keywords{Structured populations;\, diffusion;\,
Wentzell-Robin boundary
condition;\, steady states;\, spectral methods;.}
\date{\today}

\begin{abstract}
\begin{sloppypar}
We introduce a nonlinear structured population model with diffusion in the state space.
Individuals are structured with respect to a continuous variable which represents a
pathogen load. The class of uninfected individuals constitutes a special compartment that carries mass, hence
the model is equipped with ge\-ne\-ra\-lized Wentzell (or dynamic) boundary conditions. Our model is intended to describe
the spread of infection of a vertically transmitted disease, for example {\it Wolbachia} in a mosquito population.
Therefore the (infinite dimensional) nonlinearity arises in the recruitment term. First we establish global existence of solutions and the
Principle of Linearised Stability for our model. Then, in our main result, we formulate simple conditions, which guarantee
the existence of non-trivial steady states of the model. Our method utilizes an operator theoretic framework combined with a
fixed point approach. Finally, in the last section we establish a sufficient condition for the local asymptotic stability of the positive steady state.
\end{sloppypar}
\end{abstract}
\maketitle

\section{Introduction of the model}

Structured population dynamics is an exciting field of research in mathematical bi\-o\-lo\-gy,
see e.g. \cite{CUS,I,MD,WEB} for a collection of classical models and results in the area.
Structured population models often assume spatial homogeneity of the population in
a given habitat and only focus on the  dynamics of the population arising from
differences between individuals with respect to some physiological structure.
In a general model, reproduction, death and growth characterize individual
behaviour, which may be affected by competition (scramble or contest), for example for available
resources or for a mating partner. One usually incorporates certain types of nonlinearities via density dependence
in the vital rates to account for these biological phenomena.

Traditionally, structured population models have been formulated as first order
nonlinear partial differential equations of hyperbolic type and therefore the analysis of these
models is often challenging from the mathematical point of view.
In the very recent paper \cite{H} Hadeler introduced size-structured population
models with diffusion in the size-space. The biological motivation is that
diffusion allows for ``stochastic noise'' to be incorporated in the model
equations in a deterministic fashion. In \cite{H} Hadeler mainly addressed
the question that what type of boundary conditions are necessary to be imposed for a biologically
plausible and mathematically sound model. In this context particular cases
of a Robin boundary condition were considered. We note that other researchers also investigated recently 
the effects of introducing stochasticity (for example variation in individual growth rates) 
via diffusion in structured population models, see \cite{Banks1,Banks2, CDMR}.

In this paper we introduce and analyze a structured population model,
with so called distributed recruitment term and generalized Wentzell boundary conditions.
Our model intends to describe the dynamics of a population, which is infected with
a certain type of bacteria, see later for more details. At the same time our model
is general enough to make the forthcoming analysis to be applicable for
similar classes of models. In particular, we consider the following model:
\begin{align}
& u_t(x,t)+
\left(\gamma(x)u(x,t)\right)_x=\left(d(x)u_x(x,t)\right)_x-\mu(x)u(x,t) \nonumber \\
& \hspace{35mm}+\beta_0\left(U(t)\right)\int_0^m\beta_1(x,y)\beta_2\left(\frac{\int_y^mu(r,t)\,\ud r}{U(t)}\right)u(y,t)\,\ud y,\,\,
x\in(0,m),\label{equation1} \\
& \left[\left(d(x)u_{x}(x,t)\right)_x\right]_{x=0}-(d(0)+\gamma(0))u_x(0,t) \nonumber \\
& \quad+\left(\gamma(0)-\gamma'(0)-b_0\left(U(t)\right)b_2\left(\frac{\int_0^mu(r,t)\,\ud r}{U(t)}\right)\right)u(0,t)=0,
\label{equation2} \\
& d(m)u_x(m,t)-\gamma(m)u(m,t)=0,\label{equation3}
\end{align}
with a suitable initial condition.
The function $u=u(x,t)$ denotes the density of individuals of pathogen load $x$ at time $t$.
This means that the total number of individuals with pathogen (or bacterium) load between $x_1$ and $x_2$ is given by $\int_{x_1}^{x_2} u(y,t)\,\ud y$.
It is then clear that the uninfected individuals constitute a separate special class.
Therefore, in contrast to any other state, the state $x=0$ carries mass.
Hence we specify the value of the function $u$ at $x=0$, which gives the total number of individuals in the uninfected class.
Then the total population size at any time $t$ is given by
\begin{equation*}
U(t)=\left|u(0,t)\right|+\int_0^m\left|u(r,t)\right|\,\ud r=||{\bf u}||_\mathcal{X}.
\end{equation*}
In mathematical terms, we distinguish between equivalence classes of the state space $L^1$, which are represented by functions,
which have different values at $x=0$. Or, in other words, to each equivalence class of the Lebesgue space $L^1$
we assign a real number, which is its "value" at $x=0$.
This is expressed in the Wentzell-type boundary condition \eqref{equation2}.
We also assume a maximal infection load $m$. The function $\gamma$ represents the reproduction rate of the bacterium in the host.
Since the maximum infection load is $m$ we have the boundary condition \eqref{equation3},
i.e.  the infection load cannot increase above $m$ in any individual in the population.
It is shown by a straightforward calculation (see e.g. \cite{FHin3} for a similar calculation) that the boundary conditions
\eqref{equation2}-\eqref{equation3} guarantee conservation of the total population size
in the absence of mortality and recruitment. As usual, $\mu$
denotes the mortality rate (both natural mortality and extra mortality due to the infection) of individuals,
while $d$ stands for diffusion. Now we impose the following general assumptions on the model ingredients:
\begin{align*}
& \beta_0,b_0\in C([0,\infty)),\beta_2,b_2\in C([-1,1]),\mu\in C([0,m]),\quad \beta_1\in C([0,m]\times [0,m]),   \\
& \beta_0,\beta_1,\beta_2,b_0,b_2,\mu \geq 0,\quad \gamma, d\in
C^1([0,m]),\quad d>0.
\end{align*}
In addition we will make some technical assumptions on some of the
model ingredients later on. Our goal is to describe a model of a
population, which is infected with a disease that is only
transmitted vertically; hence the only nonlinearity arises in the
term specifying the recruitment of individuals. To
be more precise, we are interested in the dynamics of {\it
Wolbachia} infections in mosquito populations. {\it Wolbachia} is
a maternally transmitted endosymbiont bacterium, which infects
around 70\% of the arthropod species. Recent laboratory results
suggest, see \cite{McMeniman}, that a stable introduction of a
particular strain of {\it Wolbachia} into {\it A. Aegypti} halved
the average life span of an adult mosquito. This represents a
potential tool to eradicate mosquito born diseases (like malaria,
dengue or West Nile virus), since the life span of the mosquito
is reduced to 2-3 weeks which more or less equals the intrinsic
excubation period of the disease, i.e. the time needed for a
mosquito to become infectious after feeding on an infected human.
According to \cite{McMeniman}, high maternal inheritance, strong
cytoplasmic incompatibility and low costs to reproductive output
are the key factors for prevalence of the {\it Wolbachia}
infection. Cytoplasmic incompatibility (CI for short) means that
from a mating of an uninfected female and infected male there is
no (complete CI) or only a very few (partial CI) number of viable
offspring. We refer the interested reader to
\cite{Engelstatter,FHin2,Sinkins} for more biological background
and mathematical modelling of populations infected with a
cytoplasmic incompatibility inducing {\it Wolbachia}. The possible
underlying mechanisms for complete or partial expression of CI are
still a matter of debate for different species. In our model we
take the following view of partial CI,  see
\cite{Bourtzis,Guillemaud} for further reference. We make the
assumption that a female who has bacterium load $y$ can
successfully mate only with males of bacterium load less than $y$.
Therefore the function $\beta_2$, which determines the fertility
rate of an individual with infection load $y$, depends on the
proportion of the individuals, which have infection load higher
than $y$. In particular $\beta_2$ is a monotone decreasing
function of its argument.  In other words, we assume that the sex
ratio is $1:1$ and the probability for a female mosquito of
producing a viable offspring is a monotone decreasing function of
the proportion of the male population size that have higher
infection load. The function $\beta_0$ represents competition
effects due for example to limitations in available resources.
Similarly, $b_2$ measures the reproductive success of uninfected
females and $b_0$ represents the effects of competition on the
reproduction of uninfected individuals. The maternal transmission
rate is determined by the function $\beta_1$, i.e. infected
mothers with bacterium load $y$ give birth to offspring with
bacterium load $x$ at a rate of $\beta_1(x,y)$.

The boundary condition \eqref{equation2} can be rewritten, and
more easily understood from the point of view of its biological
meaning, in the dynamic form (by formally replacing the diffusion
operator from equation \eqref{equation1} on the boundary) as:
\begin{align}
u_t(0,t) & = u(0,t)\left(-\gamma(0)-\mu(0)+b_0(U(t))b_2\left(\frac{\int_0^mu(r,t)\,\ud r}{U(t)}\right)\right) \nonumber \\
& +u_x(0,t)d(0)+\beta_0(U(t))\int_0^m\beta_1(0,y)\beta_2\left(\frac{\int_y^mu(r,t)\,\ud r}{U(t)}\right)u(y,t)\,\ud y.\label{dynamic1}
\end{align}
It is not difficult to verify that in the absence of mortality and
recruitment, i.e. taking $\mu\equiv0$ and $b_0\equiv\beta_0\equiv
0$, the total population size $U(t)$ is conserved. Notice that in
this case, $u_t(0,t)$ equals $-\gamma(0)u(0,t)+d(0)u_x(0,t)$, i.e.,
the flux through $x=0$ due to reproduction of the bacteria (from the biological 
point of view it would be natural to assume
$\gamma(0)=0$, but we do not impose such a restriction here) 
and to diffusion (noise) in the second term (this could take into account a certain -very weak- horizontal transmission).

Our main goal in this paper is to establish sufficient conditions
for the existence of positive steady state solutions of model
\eqref{equation1}-\eqref{equation3}. We shall refer the interested
reader to \cite{CS,CS2,FGH,FHin3} where different size-structured models
with distributed recruitment processes were investigated. The
boundary condition \eqref{equation2} is the so called generalized
Wentzell-Robin or dynamic boundary condition. These ``unusual''
boundary conditions were investigated recently for models
describing physical processes such as diffusion and wave
propagation, see e.g.~\cite{FGGR,FGGR2,GG}. Briefly, they are used
to model processes where particles reaching the boundary of a
domain can be either reflected from the boundary or they can be
absorbed and then released after some time. Recently in
\cite{FHin3} we introduced, as far as we know for the first time,
Wentzell-type boundary conditions in the context of
phy\-si\-o\-lo\-gi\-cal\-ly structured populations with a
distributed recruitment process and with diffusion in the size
space. The introduction of diffusion in the size space is very
natural from the application point of view, since in the real
world individuals who start their life in the same cohort do not
finish their life so, due for example to stochastic variations in
individual growth rates. We refer the reader to \cite{H} (see also
\cite{CDMR}) where different types of population models were
introduced with diffusion in the size-space.

The idea of considering population models where the structuring variable represents a pathogen load
is clearly not new. In fact in \cite{WHG} Waldst\"{a}dter et al. introduced a similar model, where they have derived a
dynamic boundary condition for the special class of uninfected individuals. In their model this class however
does not correspond to the value $u(0,t)$, in fact the uninfected population size is represented by a new variable
$U(t)$ and a transition condition between the compartment $U$ and $u(0,t)$ is given. This implies two things:
firstly $u(0,t)\ne U(t)$, secondly individuals who are sitting in the special compartment $U$ are not subject to
diffusion, i.e. first they need to enter the state $u(0,t)$. This in some sense may seem counterintuitive, as diffusion 
is incorporated to model stochastic noise, which could result in low probability random infections. 
Since we model an (almost) completely vertically transmitted disease, we do not take into account infection in a usual way. 
By usual way we mean by means of a nonlinear transition rate of the type $SI$, which in our model would correspond to 
$u(0,t)\int_0^m u(x,t)\,\ud x$. 
But, due to the nature of the Wentzell boundary condition our model allows for transition between the uninfected and infected compartments as a diffusive flux through $x = 0$ in both directions (one of them corresponding to random infections and the other one to spontaneous, i.e. random 
healings), depending on the sign of $u_x$, as commented above.
We note that the mathematical analysis presented in \cite{WHG} is carried out in the Hilbert space
$L^2$.

We shall mention here  that the first papers introducing boundary conditions that involve second order
derivatives for parabolic or elliptic differential operators go back to the
1950s, see  the papers by Feller \cite{F1,F2} and Wentzell (also transliterated
as ``Ventcel'\! '', {\cyr Aleksandr D.~Ventcel\cyrsftsn}) \cite{V1,V2}.
These first studies were purely motivated from the mathematical point of view.
The original question, as far as we know, was to identify the maximal set of possible boundary
conditions that give rise for a parabolic differential operator to generate a
contraction semigroup on an appropriate state space.

In this paper first we establish global existence and positivity of solutions of model
\eqref{equation1}-\eqref{equation3}. Our existence proof relies on the existence of the
semigroup governing the linear part of the model, which was established in \cite{FHin3}
(with two point Wentzell boundary conditions) following similar
arguments developed in \cite{FGGR,FGGR2}, and on the analyticity of this linear semigroup, which was very recently 
established in \cite{FGGR3} (and in fact for higher dimensional domains, as well). Then in Section \ref{section:steadystates} we investigate
the existence of non-trivial steady states to our model. To treat the steady state problem,
we devise Schauder's fixed point theorem combined with an operator theoretic approach.
A similar operator theoretic framework was previously utilised for simpler problems 
(in particular with finite dimensional nonlinearities and classical boundary conditions), see e.g. \cite{CS,FGH}.
The key idea to treat the steady state problem
is to define a linear operator for a fixed environment (nonlinearity) and to study spectral properties of that operator.
Finally in the last section we establish some sufficient stability conditions for the positive steady state.

\section{Existence and positivity of solutions}\label{section:existence}

In this section we establish global existence (and positivity) of solutions of the non\-li\-ne\-ar problem \eqref{equation1}-\eqref{equation3}.
Throughout the section we employ some standard results from \cite{Henry}.
We introduce the state space $\mathcal{X}=L^1(0,m)\oplus\mathbb{R}$ with norm
\begin{equation*}
||{\bf u}||=||(u,u_0)||=||u||_1+|u_0|,
\end{equation*}
which is a Banach lattice. Next we write our problem \eqref{equation1}-\eqref{equation3} in the form of an abstract Cauchy problem:
\begin{equation}
\frac{\ud {\bf u}}{\ud t}-A{\bf u}=F({\bf u}), \quad t>0,\quad {\bf u}(0)={\bf u}_0,\label{nonlinear}
\end{equation}
where $A$ is the linear operator defined by
\begin{equation}\label{Adef}
A{\bf u}=
\begin{pmatrix}
\mathcal{A}u  \\
-\mu(0)u(0)+d(0)u'(0)-\gamma(0)u(0) \\
\end{pmatrix}.
\end{equation}
In \eqref{Adef} the operator $\mathcal{A}$ is defined as
\begin{equation}
\mathcal{A}u=\frac{\partial}{\partial x}\left(d(\cdot)\frac{\partial u}{\partial x}\right)-\frac{\partial}{\partial x}\left(\gamma(\cdot)u\right)-\mu(\cdot)u.
\end{equation}
The domain of the linear part $A$ is given by
\begin{align*}
D(A)= & \left\{ u\in C^2[0,m] \: :\: \mathcal{A} u \in L^1(0,m), \quad d(m)u'(m)-\gamma(m)u(m) =0 \right. \nonumber \\
& \left. \,\, (d(s)u'(s))'\big|_{s=0}-(d(0)+\gamma(0))u'(0)+(\gamma(0)-\gamma'(0))u(0)=0\right\}.
\end{align*}
The nonlinear but bounded function $F$ is defined as
\begin{equation}\label{Fdef}
F({\bf u})=\begin{Bmatrix}
\,\,\,F^0({\bf u})\quad \text{if}\quad {\bf u}\in\mathcal{X}\setminus\{ {\bf 0}\} \\
 {\bf 0} \quad \text{if} \quad {\bf u}={\bf 0}
\end{Bmatrix},
\end{equation}
with $D(F)=\mathcal{X}$, where
\begin{equation*}
F^0({\bf u})=\begin{pmatrix} \beta_0\left(||{\bf u}||\right)\int_0^m\beta_1(\cdot,y)\beta_2\left(\frac{\int_y^mu(r)\,\ud r}{{\bf ||u||}}\right) u(y)\,\ud y \\
\beta_0\left(||{\bf u}||\right)\int_0^m\beta_1(0,y)\beta_2\left(\frac{\int_y^mu(r)\,\ud r}{||{\bf u}||}\right) u(y)\,\ud y +b_0\left(||{\bf u}||\right)b_2\left(\frac{\int_0^mu(y)\,\ud y}{||{\bf u}||}\right)u_0\\
\end{pmatrix}.
\end{equation*}
Note that $F^0$ is only defined on the open set $\mathcal{X}\setminus\{{\bf 0}\}$.

\begin{theorem}\label{existencesol}
Assume that $\beta_0,\beta_2,b_0$ and $b_2$ are locally
Lipschitzian. Then a unique solution of problem \eqref{nonlinear}
exists for all positive times for any $u_0 \in \mathcal{X}$.
\end{theorem}
\noindent {\bf Proof.}
We use \cite[Corollary 3.3.5]{Henry} to establish global existence of solutions to the semilinear problem \eqref{nonlinear}.
To this end, first we show that $F$ is locally Lipschitz continuous, i.e. for every $\mathbf{u}\in\mathcal{X}$ there exists a neighbourhood
$U\subset\mathcal{X}$ of $\mathbf{u}$  such that for every
${\bf u}_1,{\bf u}_2\in U$ we have
$\left|\left|F({\bf u}_1)-F({\bf u}_2)\right|\right|_\mathcal{X}\le L\,\left|\left|{\bf u}_1-{\bf u}_2\right|\right|_\mathcal{X}
=L\left(||u_1-u_2||_1+|u_{1_0}-u_{2_0}|\right)$ for some constant $L\ge 0$.

If either ${\bf u}_1={\bf 0}$ or ${\bf u}_2={\bf 0}$
then since $F$ is bounded (i.e. all the ingredients $\beta_0,\beta_1,\beta_2,b_0$ and $b_2$ are continuous, hence bounded) and $F({\bf 0})={\bf 0}$
we have the Lipschitz property $||F({\bf u})||\le L||{\bf u}||$ for some $L>0$. For ${\bf u}_1\ne{\bf 0},{\bf u}_2\ne {\bf 0}$ we have:
\begin{align}
&  \left|\left|F({\bf u}_1)-F({\bf u}_2)\right|\right|_\mathcal{X}  \le \nonumber \\
& \quad \left|\beta_0\left(||{\bf u}_1||\right)-\beta_0\left(||{\bf u}_2||\right)\right| \left|\left|\int_0^m\beta_1(\cdot,y)\beta_2\left(\frac{\int_y^m u_1(r)\,\ud r}{||{\bf u}_1||}\right)u_1(y)\,\ud y \right|\right|_1\label{Leq1} \\
& \quad +\beta_0\left(||{\bf u}_2||\right)\left|\left|\int_0^m\beta_1(\cdot,y)\left(\beta_2\left(\frac{\int_y^mu_1(r)\,\ud r}{||{\bf u}_1||}\right)u_1(y)-\beta_2\left(\frac{\int_y^mu_2(r)\,\ud r}{||{\bf u}_2||}\right)u_2(y)\right)\,\ud y\right|\right|_1 \label{Leq2} \\
& \quad + \left|\beta_0\left(||{\bf u}_1||\right)-\beta_0\left(||{\bf u}_2||\right)\right| \left|\int_0^m\beta_1(0,y)\beta_2\left(\frac{\int_y^m u_1(r)\,\ud r}{||{\bf u}_1||}\right)u_1(y)\,\ud y \right| \label{Leq3} \\
& \quad +\beta_0\left(||{\bf u}_2||\right)\left|\int_0^m\beta_1(0,y)\left(\beta_2\left(\frac{\int_y^mu_1(r)\,\ud r}{||{\bf u}_1||}\right)u_1(y)-\beta_2\left(\frac{\int_y^mu_2(r)\,\ud r}{||{\bf u}_2||}\right)u_2(y)\right)\,\ud y\right| \label{Leq4} \\
& \quad +\left|b_0\left(||{\bf u}_1||\right)-b_0\left(||{\bf u}_2||\right)\right|b_2\left(\frac{\int_0^m u_1(r)\,\ud r}{||{\bf u}_1||}\right)u_{1_0} \label{Leq5} \\
& \quad +b_0\left(||{\bf u}_2||\right)\left(b_2\left(\frac{\int_0^m u_1(r)\,\ud r}{||{\bf u}_1||}\right)u_{1_0}-b_2\left(\frac{\int_0^m u_2(r)\,\ud r}{||{\bf u}_2||}\right)u_{2_0}\right).\label{Leq6}
\end{align}
The term in \eqref{Leq1} can be bounded above by
$L_{\beta_0}||{\bf u}_1-{\bf
u}_2||M_{\beta_1}M_{\beta_2}m||u_1||_1$, where
$\beta_0<M_{\beta_0}, \beta_1<M_{\beta_1},\beta_2<M_{\beta_2}$ and
$L_{\beta_0}$ is the Lipschitz constant of $\beta_0$ in
$\left[0,\displaystyle\sup_{u\in U}||{u}||\right]$. To obtain the estimate for
\eqref{Leq2} we note that
\begin{align}
&\left| \beta_2\left(\frac{\int_y^mu_1(r)\,\ud r}{||{\bf u}_1||}\right)u_1(y)-\beta_2\left(\frac{\int_y^mu_2(r)\,\ud r}{||{\bf u}_2||}\right)u_2(y)\right| \nonumber \\
& =\left|\left(\beta_2\left(\frac{\int_y^mu_1(r)\,\ud r}{||{\bf u}_1||}\right)-\beta_2\left(\frac{\int_y^mu_2(r)\,\ud r}{||{\bf u}_2||}\right)\right)u_1(y)+\beta_2\left(\frac{\int_y^mu_2(r)\,\ud r}{||{\bf u}_2||}\right)(u_1(y)-u_2(y))\right| \nonumber \\
& \le L_{\beta_2}\left|\left|\frac{\int_y^mu_1(r)\,\ud r}{||{\bf u}_1||}-\frac{\int_y^mu_2(r)\,\ud r}{||{\bf u}_1||}+\frac{\int_y^mu_2(r)\,\ud r}{||{\bf u}_1||}-\frac{\int_y^mu_2(r)\,\ud r}{||{\bf u}_2||}\right|\right|_1\left|u_1(y)\right|+M_{\beta_2}|u_1(y)-u_2(y)|\nonumber \\
& \le L_{\beta_2}\left(\left|\left|\int_y^mu_1(r)\,\ud r-\int_y^mu_2(r)\,\ud r\right|\right|_1\frac{\left|u_1(y)\right|}{||{\bf u}_1||}\right)+M_{\beta_2}|u_1(y)-u_2(y)|\nonumber \\
& \quad+L_{\beta_2}\left(\Big| ||{\bf u}_2||-||{\bf
u}_1||\Big|\frac{\left|\left|\int_y^mu_2(r)\,\ud
r\right|\right|_1}{||{\bf u}_2||}\frac{\left|u_1(y)\right|}{||{\bf
u}_1||}\right),\label{Leq7}
\end{align}
where $L_{\beta_2}$ is the Lipschitz constant of $\beta_2$ in $[-1,1]$. Similarly to \eqref{Leq7}, we obtain the appropriate estimates for the terms
\eqref{Leq3}-\eqref{Leq6}. Hence $F$ is locally Lipschitzian. Since $||F({\bf u})-F({\bf 0})||=||F({\bf u})||\le L||{\bf u}||$, $F$
is also sublinear and since $A$ is sectorial (see \cite{FGGR3}) the statement of Theorem \ref{existencesol} follows from \cite[Corollary 3.3.5]{Henry}.
\hfill $\Box$

Note that $F\,:\, \mathcal{X}_+\to\mathcal{X}_+$. The solution of the Cauchy problem \eqref{nonlinear} can be written as
\begin{equation}\label{nonlinear2}
{\bf u}(t)=\mathcal{T}(t){\bf u}_0+\int_0^t\mathcal{T}(t-s)F({\bf u}(s))\,\ud s,\quad t\in (0,t_0),
\end{equation}
where $\mathcal{T}$ is the linear semigroup generated by the closure of the sectorial operator $A$. We refer the reader to
\cite{FHin3} were it is showed that a similar linear problem (with two point Wentzell boundary condition) is governed by a quasicontractive positive semigroup. Furthermore we refer to \cite{FGGR3}  where the sectoriality of the operator $A$ was shown in $L^p,\, 1\le p\le \infty$, in general. Since $F$ is a positive operator and the semigroup $\mathcal{T}(t)$ is positive the variation formula \eqref{nonlinear2} immediately shows positivity of solutions.

The Principle of Linearised Stability can be established for model \eqref{equation1}-\eqref{equation3} by directly applying again results from \cite{Henry}.
In particular, if ${\bf u}_*$ is an equilibrium point then let $F'_{{\bf u}_*}$ denote the linearisation of $F$ at the equilibrium ${\bf u}_*$, 
which is a bounded li\-ne\-ar operator defined on $\mathcal{X}$ (if it exists). See later in Section 4 for more details.
Then \cite[Theorem 5.1.1]{Henry} implies that if $A$ is sectorial and $F$ is locally Lipschitz continuous then
the equilibrium ${\bf u}_*$ is asymptotically stable if
$\sigma\left(A+F'_{{\bf u}_*}\right)\subset\left\{\lambda\in\mathbb{C}\,|\, Re(\lambda)<\alpha<0\right\}$.
On the other hand, \cite[Theorem 5.1.3]{Henry} implies (under the same conditions on $A$ and $F$) that the equilibrium ${\bf u}_*$ is unstable if
$\sigma\left(A+F'_{{\bf u}_*}\right)\cap\left\{\lambda\in\mathbb{C}\,|\, Re(\lambda)>0\right\}$ is not empty.

\section{Existence of non-trivial steady states}\label{section:steadystates}

It is obvious that model \eqref{equation1}-\eqref{equation3}
admits the trivial steady state. It is also clear that even the
time independent version of equations
\eqref{equation1}-\eqref{equation3} cannot be solved explicitly.
Therefore, to establish conditions which guarantee the existence
of a positive steady state we utilise a combination of an operator
theoretic framework (see e.g. \cite{CS,FGH}) and a fixed point
approach. For basic definitions and results from linear semigroup
theory used throughout this section we refer the reader to
\cite{AGG,CH,NAG}.

For a fixed non vanishing ${\bf v}\in L^1_+(0,m)\oplus\mathbb{R}$ let us define the linear operator $\Psi_{{\bf v}}$ (parametrised by ${\bf v}$) by
\begin{equation}
\Psi_{{\bf v}}{\bf u}=A{\bf u}+F^0_{\bf v}{\bf u},\label{Psioperator}
\end{equation}
where the operator $A$ is defined in \eqref{Adef}, and
\begin{equation*}
F^0_{\bf v}{\bf u}=\begin{pmatrix} \beta_0\left(||{\bf v}||\right)\int_0^m\beta_1(\cdot,y)\beta_2\left(\frac{\int_y^m v(r)\,\ud r}{{\bf ||v||}}\right) u(y)\,\ud y \\
\beta_0\left(||{\bf v}||\right)\int_0^m\beta_1(0,y)\beta_2\left(\frac{\int_y^mv(r)\,\ud r}{||{\bf v}||}\right) u(y)\,\ud y +b_0\left(||{\bf v}||\right)b_2\left(\frac{\int_0^mv(y)\,\ud y}{||{\bf v}||}\right)u_0\\
\end{pmatrix},
\end{equation*}
with domain
\begin{align}
D(\Psi_{\bf v})= & \left\{ u\in C^2([0,m]) \: :\: \mathcal{A}u \in L^1(0,m),\quad d(m)u'(m)-\gamma(m)u(m) =0 \right. \nonumber \\
& \left. (d(s)u'(s))'\big|_{s=0}-(d(0)+\gamma(0))u'(0)+(\gamma(0)-\gamma'(0))u(0)=0\right\}. \label{Psidef}
\end{align}
It is clear that $(v_*,v_*(0))$ is a nontrivial steady state of model \eqref{equation1}-\eqref{equation3}
if (and only if) it is a positive eigenvector belonging to the kernel of the (closure of the) linear operator $\Psi_{\bf v}$.

Next we establish some (necessary) regularity properties of the
semigroup generated by the closure of $\Psi_{\bf v}$. Several
characterisations of irreducibility of a positive semigroup exist
in the literature, see for example in \cite{AGG,CH,NAG}. Here we
follow \cite{NAG}. A positive semigroup $\mathcal{T}$ on the Banach lattice $\mathcal{X}$ is said to
be irreducible if the resolvent of its generator $G$ is strictly
positive, i.e. $\forall \,0\not\equiv f\in \mathcal{X}_+$ we have
$(R(\lambda,G)f)(x)>0$ for almost all $x$ and some $\lambda>s(G)$.
 \begin{lemma}\label{posirred}
For every fixed ${\bf v}\in \mathcal{X}_+$ the closure of the linear operator $\Psi_{\bf v}$ generates an irreducible semigroup on $\mathcal{X}$.
\end{lemma}
\noindent {\bf Proof.} We introduce the mortality operator $M$ as follows:
\begin{align}
& M{\bf u}=\begin{pmatrix} -\mu(\cdot)u \\ -\mu(0)u(0) \\ \end{pmatrix} \quad \text{on}\quad \mathcal{X}_+. \label{operators}
\end{align}
We then consider the resolvent equation
\begin{equation}
\left(\lambda\mathcal{I}-(A-M)\right) {\bf u}={\bf h},
\end{equation}
for ${\bf h}\in \mathcal{X}_+,\,\lambda>0$, i.e.
\begin{align}
-h(x) & =\left(d(x)u_x(x)-\gamma(x)u(x)\right)_x-\lambda u(x),\quad x\in (0,m), \label{reseq1} \\
-h_0 & =-\left(\gamma(0)+\lambda\right)u(0)+d(0)u'(0),\label{reseq2}
\end{align}
for an unknown $u\in D(A).$ We have shown in \cite[Theorem 2.1]{FHin3} that the linear semigroup generated by $A-M$ is positive (with two point Wentzell boundary conditions, but the proof
can be trivially adapted to our case), i.e. the resolvent operator $\left(\lambda\mathcal{I}-(A-M)\right)^{-1}$ is positive for $\lambda>0$ large enough.
Hence the solution ${\bf u}$ of equations \eqref{reseq1}-\eqref{reseq2} is non-negative, i.e. ${\bf u}\ge 0$. It is only left to show
that the solution ${\bf u}$ is in fact strictly positive.

The one-dimensional minimum principle assures that $u$ cannot
attain its minimum value at an interior point of $(0,m)$, see e.g.
\cite[Ch. 1 Theorem 3]{PW}. On the other hand, \cite[Ch. 1 Theorem 4]{PW} guarantees that if $u$ attains its non-positive
minimum at $0$ then $u'(0)>0$ holds, which contradicts ${\bf h}\ge
0$, whereas the boundary condition at $x=m$ implies $u(m)>0$ by
the same theorem. Therefore the semigroup generated by the closure
of $A-M$ is irreducible. Since $M$ is a bounded multiplication
operator and $F^0_{\bf v}$ is bounded and positive for every ${\bf
v}\in\mathcal{X}_+$ it follows that the closure of $\Psi_{\bf v}$
generates an irreducible semigroup for every ${\bf
v}\in\mathcal{X}_+$, see e.g. \cite[C-III Proposition 3.3]{AGG}. \hfill $\Box$

\begin{remark}
First of all we note that every steady state $(v,v(0))$ of model
\eqref{equation1}-\eqref{equation3} shall have regularity
$W^{2,1}(0,m)$. This immediately implies that model
\eqref{equation1}-\eqref{equation3} does not admit a steady state
which has uninfected individuals only, i.e. of the form
$(0,v(0))$, with $v(0)\ne 0$.

Secondly, as we noted before, a vector ${\bf v}=(v,v(0))$  is a steady state of \eqref{equation1}-\eqref{equation3} if and only if
it is a positive eigenvector belonging to the kernel of the closure of $\Psi_{\bf v}$.
Therefore, irreducibility of the semigroup generated by the closure of  $\Psi_{\bf v}$
shows that we also cannot have a steady state of the form $(v,0)$,
since every positive eigenvector of the closure of the generator is strictly positive if the semigroup is irreducible.
\end{remark}

\begin{lemma}\label{compact}
For every non vanishing ${\bf v}\in \mathcal{X}_+$ the spectrum of $\Psi_{\bf v}$ can contain only isolated eigenvalues of finite algebraic multiplicity.
\end{lemma}
\noindent {\bf Proof.} Since $M$ and $F^0_{\bf v}$ are bounded it
is enough to show that $R(\lambda,A-M)$ is compact. This follows
however by noting that the solution of the resolvent equation
\eqref{reseq1}-\eqref{reseq2} is in $W^{2,1}(0,m)\oplus\mathbb{R}$
which is compactly embedded in $\mathcal{X}$. The statement now
follows on the grounds of \cite[Proposition II.4.25]{NAG} and
\cite[Corollary IV.1.19]{NAG}. \hfill $\Box$

As we noted before our goal is to show that there exists an element ${\bf v}\in\mathcal{X}_+$ such that the operator $\Psi_{\bf v}$ has eigenvalue zero.
Then, Lemmas \ref{posirred} and \ref{compact} guarantee the existence of a corresponding strictly positive (unique normalized) eigenvector.
In case of a model with one dimensional nonlinearity, such as the one we treated in \cite{FGH}, the operator $\Psi$ is in fact parametrized by a scalar
quantity, namely, the total population size. In this case the positive steady state is obtained readily after an appropriate normalization
of the positive eigenvector. In case of the model treated here the function $\beta_2$ naturally depends on an infinite dimensional variable.
Therefore, in general there is no guarantee that the positive eigenvector corresponding to the zero eigenvalue of the operator $\Psi_{\bf v}$
is in fact ${\bf v}$ (or a scalar multiple of it). For this reason, we need to construct an appropriate nonlinear map
on a certain level set of the positive cone of the state space $\mathcal{X}$ and establish the existence of a fixed point of this map.
We note that a different fixed point strategy was employed very recently in \cite{FHin4} for a class of models with
infinite dimensional nonlinearities (and zero flux boundary condition).
That method, which uses fixed point results of nonlinear maps in conical shells of Banach spaces, does not
apply to our model here since the construction of the nonlinear map requires the (implicit) solution of the steady state equation. We also note that a similar argument to the one used here, even though in the case
of finite dimensional interaction variable is used in \cite{Pal}.

\begin{theorem}\label{existencess}
Assume that $\beta_0$ and $b_0$ are strictly monotone decreasing functions of their argument and
\begin{enumerate}
\item[(i)]
$\displaystyle\lim_{x\to\infty}\beta_0(x)=\displaystyle\lim_{x\to\infty}b_0(x)=0$ \quad and\quad $\mu(x)>\mu_0>0,\quad \forall x\in (0,m)$; \vspace{2mm}
\item[(ii)] there exists an $r>0$ such that for all ${\bf v_*}\in\mathcal{X}_+$, $||{\bf v_*}||\le r$ we have that the
spectral bound $s\left(\Psi_{\bf v_*}\right)>0$.
\end{enumerate}
Then model \eqref{equation1}-\eqref{equation3} has at least one non-trivial steady state.
\end{theorem}
\noindent {\bf Proof.} Since $\Psi_{\bf v}$ is a generator of a
positive and irreducible semigroup with compact resolvent it
follows that the spectrum $\sigma(\Psi_{\bf v})$ is not empty, see
e.g. \cite[C-III Theorem 3.7]{AGG}. Moreover, the spectral bound
$s(\Psi_{\bf v})$ is an isolated eigenvalue of algebraic
multiplicity one with a corresponding strictly positive
eigenvector and it is the unique eigenvalue with positive
eigenvector, see e.g. \cite[Theorem 9.10]{CH}. We also note that
the spectral bound and the corresponding positive eigenvector
change continuously with respect to (the parameter) ${\bf v}$, see
e.g. \cite[Sect. 3 in Ch.4]{K}.

We introduce the level set
\begin{equation*}
S=\left\{{\bf x}\in \mathcal{X}_+\,| \, s(\Psi_{\bf x})=0\right\}.
\end{equation*}
It is shown that the closure of the linear operator $A-M$ generates a (positive) contraction semigroup, hence we have $s(A-M)=\omega_0\le 0$.
Intuitively it is clear that contractivity follows simply because in the absence of mortality and recruitment the total population size is preserved. Therefore, if
\begin{equation*}
\left|\left|F^0_{\bf v}\right|\right|<\displaystyle\inf_{x\in [0,m]}\{\mu(x)\},
\end{equation*}
then it is shown that $s\left(\Psi_{\bf v}\right)<0$. Hence
conditions (i) imply that there exists an $R>0$ such that for
$\forall\,{\bf v}_*\in\mathcal{X}_+$, $||{\bf v_*}||\ge R$ we have
$s\left(\Psi_{\bf v_*}\right)<0$. From (ii) it then follows that
$S \subset \left\{{\bf x}\in \mathcal{X}_+\,|\, r <
 ||{\bf x}|| < R \right \}.$ It also follows from conditions
(i) and (ii) and from the continuity and strict monotonicity of
the spectral bound (see below) that along every positive ray
$\mathcal{R}=\left\{\alpha {\bf v_*}\,|\, \alpha\in
\mathbb{R}_+,\, {\bf v_*}\in \mathcal{X}_+\setminus\left\{\,{\bf 0}\right\}\right\}$ 
there exists a (unique!) ${\bf v}$ such that
$s(\Psi_{\bf v})=0$. So the set $S$ intersects every positive ray
$\mathcal{R}$ in a unique element.

Next, we shall show that the spectral bound is in fact strictly
monotone decreasing along every positive ray. To this end let
$0<\alpha_1<\alpha_2$ be real numbers and let ${\bf v}\in
\mathcal{X}_+$. Consider the operators $\Psi_{\alpha_1{\bf v}}$
and $\Psi_{\alpha_2 {\bf v}}$. Both of them have compact
resolvents by Lemma \ref{compact} and generate positive and
irreducible semigroups. We also note that $\Psi_{\alpha_1 {\bf
v}}-\Psi_{\alpha_2{\bf v}}$ is a positive operator, since
\begin{align*}
\left(\Psi_{\alpha_1{\bf v}}-\Psi_{\alpha_2{\bf v}}\right){\bf u} & =\left(F^0_{\alpha_1{\bf v}}-F^0_{\alpha_2{\bf v}}\right){\bf u} \\
& =\begin{pmatrix} \left( \beta_0\left(\alpha_1 ||{\bf v}||\right)-\beta_0\left(\alpha_2 ||{\bf v}||\right) \right)\int_0^m\beta_1(\cdot,y)\beta_2\left(\frac{\int_y^m v(r)\,\ud r}{{\bf ||v||}}\right) u(y)\,\ud y \\
\left( \beta_0\left(\alpha_1 ||{\bf
v}||\right)-\beta_0\left(\alpha_2 ||{\bf v}||\right) \right)
\int_0^m\beta_1(0,y)\beta_2\left(\frac{\int_y^mv(r)\,\ud r}{||{\bf v}||}\right) u(y)\,\ud y \\
\end{pmatrix} \\
&\quad +\begin{pmatrix} 0 \\
\left(b_0\left(\alpha_1||{\bf v}||\right)-b_0\left(\alpha_2||{\bf v}||\right)\right)b_2\left(\frac{\int_0^mv(y)\,\ud y}{||{\bf v}||}\right)u_0\\
\end{pmatrix} \\
& \ge 0,
\end{align*}
for every ${\bf u}\in \mathcal{X}_+$. We also have
\begin{equation*}
R(\lambda,\Psi_{\alpha_1{\bf v}})-R(\lambda,\Psi_{\alpha_2 {\bf
v}})=R(\lambda,\Psi_{\alpha_1 {\bf v}}) \left(\Psi_{\alpha_1{\bf
v}}-\Psi_{\alpha_2{\bf v}}\right)R(\lambda,\Psi_{\alpha_2 {\bf
v}})\ge 0,
\end{equation*}
for $\lambda$ large enough. Hence \cite[Proposition A.2]{AB}
implies that $s(\Psi_{\alpha_1{\bf v}})>s(\Psi_{\alpha_2{\bf
v}})$.

Next we construct the nonlinear map $\Phi\,:\, S \to S$, $\Phi({\bf v})=(f_4\circ f_3\circ f_2\circ f_1)({\bf v})={\bf v'}$,
as illustrated briefly in the following diagram:
\begin{equation}
\Phi\,:\,\underbrace{{\bf v}}_{\in S}\hspace{6mm}\xrightarrow{f_1}\underbrace{\Psi_{\bf v}}_{\in L(D(\Psi_{\bf v}),\mathcal{X})}
\xrightarrow{f_2}\underbrace{\overline{{\bf V}}_{\bf v}}_{\in W^{2,1}_+(0,m)\oplus\mathbb{R}}
\xrightarrow{f_3}\underbrace{\overline{{\bf V}}_{\bf v}}_{\in L^1_+(0,m)\oplus\mathbb{R}}
\xrightarrow{f_4}\underbrace{{\bf v'}}_{\in S}.
\end{equation}
For every element ${\bf v}\in S$ the map $f_1$ assigns the
corresponding linear operator $\Psi_{\bf v}$ "parametrized" by
${\bf v}$, which is defined via \eqref{Psioperator}. This map
$f_1$ is continuous and bounded. The map $f_2$ assigns the
strictly positive normalized eigenvector corresponding to the zero
eigenvalue of the operator $\Psi_{\bf v}$. This map $f_2$ is
clearly bounded and it is also continuous, in fact it is even
analytic, see \cite[Lemma 1.3]{CR}. The map $f_3$ is the compact
injection of $W^{2,1}(0,m)\oplus\mathbb{R}$ into
$L^1(0,m)\oplus\mathbb{R}$. Finally $f_4$ is the projection along
positive rays of the eigenvector $\overline{\bf V}_{\bf v}$ back
into the set $S$, which is again continuous and bounded.

Next we apply Schauder's fixed point theorem to an appropriately defined map.
As we noted above the map $\Phi$ is continuous and compact. It is only left to show that $S$ is homeomorphic to a convex set.
To this end, we define the map  $h\,:\, S\to B$ by $h({\bf u})=\frac{{\bf u}}{||{\bf u}||}$, where $B$ is the unit sphere intersected
with the positive cone $\mathcal{X}_+$, i.e.
\begin{equation*}
B:=\left\{{\bf u}=(u,u^0)\in\mathcal{X}_+\,\vert\, ||{\bf u}||=||u||_1+|u^0|=1\right\}.
\end{equation*}
Then $h$ is clearly continuous and one to one, since the spectral bound $s\left(\Psi\right)$ is strictly monotone decreasing along positive rays in $\mathcal{X}$
and every positive ray of $\mathcal{X}$ intersects $B$ in a unique element.

We shall show now that the function $h^{-1}\,:\, B\to S$ defined
via $h^{-1}({\bf w})=\alpha\,{\bf w}$, $\alpha\in\mathbb{R}_+$
such that $s\left(\Psi_{\alpha\, {\bf w}}\right)=0$,  is also
continuous. Let ${\bf w}_n\in B$ be a sequence such that ${\bf
w}_n\to {\bf w}\in B$ and consider the sequence $h^{-1}({\bf
w}_n)=\alpha_n{\bf w}_n$. It follows from condition (i) that there
exists an $R>0$ such that $s(\Psi_{\bf v})<0$ for every $||{\bf
v}||>R$, hence the sequence $\alpha_n$ is bounded. Let
$\alpha_{n_k}$ be a convergent subsequence of $\alpha_n$ and let
$\alpha_{n_k}\to \bar{\alpha}$. Since
$s\left(\Psi_{\alpha_{n_k}{\bf w}_{n_k}}\right)=0$ for every
$k\in\mathbb{N}$ it follows from the continuity of the spectral
bound that $s\left(\Psi_{\bar{\alpha}{\bf w}}\right)=0$. Since
there is exactly one element on the positive ray spanned by ${\bf
w}$ at which the spectral bound vanishes we have
$\bar{\alpha}=\alpha$. If $\alpha_{n_l}$ is another convergent
subsequence of $\alpha_n$ then the continuity of the spectral
bound implies again that $\alpha_{n_l}\to\alpha$. So
$\alpha_{n}\to\alpha$ and finally $h^{-1}({\bf w}_n) =
\alpha_n{\bf w}_n \to \alpha {\bf w} = h^{-1}({\bf w})$.

We note that the set $B$ is convex since we are in an AL-space, i.e. we have $||{\bf f}+{\bf g}||=||{\bf f}||+||{\bf g}||$ for every ${\bf f},{\bf g}\in\mathcal{X}_+$,
see e.g. \cite{AGG}. Finally we apply Schauder's fixed point theorem (see e.g. \cite{GT}) to the continuous and compact map
$\bar{\Phi}\,:\,B\to B$ defined by
\begin{equation*}
\bar{\Phi}({\bf x})=h\circ\Phi\circ h^{-1} ({\bf x}),
\end{equation*}
to obtain a fixed point $\bar{{\bf x}}_*\in B$ of the map $\bar{\Phi}$, which yields a fixed point ${\bf x}_*=h^{-1}\left(\bar{{\bf x}}_*\right)$ of
the map $\Phi$ in $S$. This ${\bf x}_*$ is the positive steady state of model \eqref{equation1}-\eqref{equation3} (see the characterization of a nontrivial steady state immediately after equation \eqref{Psidef}).    \hfill $\Box$

\begin{remark}
Condition (i) is natural from the biological point of view. The first condition in (i) requires that the fertility rate of both infected and uninfected individuals
tends to zero as the population size tends to infinity. This may be due for example to competition effects. In fact it turns out that as expected, cytoplasmic incompatibility
itself does not have a negative feedback on population growth. The assumption of a strictly positive mortality function in (i) seems also realistic.

Condition (ii) seems to be also a natural one, if one is to expect the existence of a positive steady state.
In fact, if it is not satisfied then one can show that there exists a monotone decreasing sequence of positive real numbers $r_n\to 0$,
such that for every $n\in\mathbb{N}$ there exists a ${\bf u}_*^n$ with $||{\bf u}_*^n||\le r_n$ and $s\left(\Psi_{{\bf u}_*^n}\right)\le 0$ holds.
Then, the continuity of the spectral bound implies that $s\left(\Psi_{\bf 0}\right)\le 0$ and since the spectral bound is a strictly monotone decreasing function along
positive rays, we have $s\left(\Psi_{\bf v}\right)\le 0$, for every ${\bf v}\in\mathcal{X}_+$.
\end{remark}

\section{Stability}
In the previous section we established conditions which guarantee the existence of a non-trivial steady state. 
The next natural step is to study the stability of  the steady state (and also the stability of the trivial steady state).  In this section we are going to establish a sufficient condition for the local asymptotic stability of the non-trivial steady state.  First we note that the operator $F$ is not (Fr\'{e}chet) differentiable at ${\bf 0}$ because it is
homogeneous of degree $1$ (it is an "angular operator"). Nevertheless it is G\^{a}teaux-differentiable at ${\bf 0}$. In particular its directional derivative at ${\bf 0}$ into the direction of any element ${\bf v}$ in the positive cone is:
\begin{align*}
& \ud F_{{\bf 0}}{\bf v}= \begin{pmatrix}  \beta_0(0)\int_0^m\beta_1(\cdot,y)\beta_2\left(\frac{\int_y^m v(r)\,\ud r}{||{\bf v}||}\right)v(y)\,\ud y \\  
\beta_0(0)\int_0^m\beta_1(0,y)\beta_2\left(\frac{\int_y^m v(r)\,\ud r}{||{\bf v}||}\right)v(y)\,\ud y+b_0(0)b_2\left(\frac{\int_0^mv(y)\,\ud y}{||{\bf v}||}\right)v_0
\end{pmatrix}, 
\end{align*}
which is a nonlinear operator. 

If ${\bf u}_*$ is a non-trivial equilibrium of \eqref{equation1}-\eqref{equation3}  then we can formally linearise equation \eqref{equation1} around ${\bf u}_*$. In particular if we denote by $F'_{{\bf u}_*}$ 
the linearisation (the Fr\'{e}chet de\-ri\-va\-tive) of the nonlinear operator $F$ at ${\bf u}_*$ then the linearised problem can be cast in the form of an abstract Cauchy problem
\begin{equation}\label{linearised}
\frac{\ud}{\ud t}\,{\bf v}=\left(A+F'_{{\bf u}_*}\right){\bf v},\quad {\bf v}(0)={\bf v}_0,
\end{equation}
where $A$ is defined in \eqref{Adef} and $D(F'_{{\bf u}_*})=\mathcal{X}$. We leave it for future work to 
address the linear problem \eqref{linearised} (i.e. to address stability questions of equilibria), in general, for example using the Liapunov function 
techniques elaborated in \cite[Sect. 4.2]{WEB}. 
Here however, we establish a straightforward and simple condition which guarantees that the non-trivial equilibrium is locally asymptotically stable. 

To this end we consider on the state space $\mathcal{X}$ the operator $F$ defined in \eqref{Fdef} with domain $D(F)=\mathcal{X}^+$. 
Then using approximations (in the spirit of \cite[Def.2.4 in Sect. 2.6]{WEB}) such as 
\begin{equation*}
\beta_0(||{\bf u}_*+{\bf v}||)\sim\beta_0(||{\bf u}_*||)+\beta'_0(||{\bf u}_*||)\left(\int_0^mv(y)\,\ud y+v_0\right)
\end{equation*}
and similarly
\begin{align*}
 \beta_2\left(\frac{\int_y^mu_*(r)+v(r)\,dr}{||{\bf u}_*+{\bf v}||}\right) & \sim\beta_2\left(\frac{\int_y^mu_*(r)\,dr}{||{\bf u}_*||}\right)   \\ & +\beta'_2\left(\frac{\int_y^mu_*(r)\,dr}{||{\bf u}_*||}\right)\left(\frac{\int_y^mv(r)\,dr}{||{\bf u}_*||}-\int_y^mu_*(r)\,dr\frac{\int_0^mv(y)\,dy+v_0}{||{\bf u}_*||^2}\right),
\end{align*}
it is shown that the linearisation of the operator $F$ at the non-trivial equilibrium ${\bf u}_*$ (according again to \cite[Def. 2.4 in Sect. 2.6]{WEB} is
\begin{align*}
& F'_{{\bf u}_*}{\bf v}= \\
 & \begin{pmatrix}  \int_0^m\beta_1(\cdot,y)\beta_2\left(\frac{\int_y^mu_*(r)\,dr}{||{\bf u}_*||}\right)\left(\beta_0(||{\bf u}_*||)v(y)+\beta'_0(||{\bf u}_*||)\left(\int_0^m v(y)\, dy+v_0\right)u_*(y)\right)\,\ud y \\  
\int_0^m\beta_1(0,y)\beta_2\left(\frac{\int_y^mu_*(r)\,dr}{||{\bf u}_*||}\right)\left(\beta_0(||{\bf u}_*||)v(y)+\beta'_0(||{\bf u}_*||)\left(\int_0^m v(y)\, dy+v_0\right)u_*(y)\right)\,\ud y  \\ 
\end{pmatrix}  \\
+ & \begin{pmatrix} \int_0^m\beta_1(\cdot,y)\beta_0(||{\bf u}_*||)\beta'_2\left(\frac{\int_y^mu_*(r)\,dr}{||{\bf u}_*||}\right)\left(\frac{\int_y^mv(r)\,dr}{||{\bf u}_*||}-\int_y^mu_*(r)\,dr\frac{\int_0^mv(y)\,dy+v_0}{||{\bf u}_*||^2}\right)u_*(y)\,\ud y \\  
 \int_0^m\beta_1(0,y)\beta_0(||{\bf u}_*||)\beta'_2\left(\frac{\int_y^mu_*(r)\,dr}{||{\bf u}_*||}\right)\left(\frac{\int_y^mv(r)\,dr}{||{\bf u}_*||}-\int_y^mu_*(r)\,dr\frac{\int_0^mv(y)\,dy+v_0}{||{\bf u}_*||^2}\right)u_*(y)\,\ud y\\ 
\end{pmatrix} \\
+ & \begin{pmatrix} 0\\  
 b_2\left(\frac{\int_0^mu_*(y)\,dy}{||{\bf u}_*||}\right)\left[b_0(||{\bf u}_*||)v_0+b'_0(||{\bf u}_*||)\left(\int_0^mv(y)\, dy+v_0\right)u_{*_{0}}\right] \\ 
\end{pmatrix} \\
+ & \begin{pmatrix} 0\\  
b_0(||{\bf u}_*||)b'_2\left(\frac{\int_0^mu_*(y)\,dy}{||{\bf u}_*||}\right)\left(\frac{\int_0^mv(y)\, dy}{||{\bf u}_*||}-\frac{\int_0^mu_*(y)\, dy}{||{\bf u}_*||^2}\left(\int_0^m v(y)\, dy +v_0\right)\right)u_{*_{0}} \\ 
\end{pmatrix}, 
\end{align*}
which is a bounded linear operator on $\mathcal{X}$.
\begin{theorem}\label{stability}
If  
\begin{equation}
\nu:=\displaystyle\inf_{s\in [0,m]}\left\{\mu(s)\right\}>\left|\left|F'_{{\bf u}_*}\right|\right| ,
\end{equation}
then the non-trivial steady state ${\bf u}_*$ is locally asymptotically stable.
\end{theorem}
\noindent {\bf Proof.} 
Note that $A-M$ has compact resolvent (where $M$ is the mortality operator introduced in \eqref{operators}), see the proof of Lemma \ref{compact}, 
and generates a positive and irreducible semigroup, see the proof of Lemma \ref{posirred}. Therefore its point spectrum $\sigma_P(A-M)$ is not empty, 
see \cite[C-III Theorem 3.7]{AGG}. We also noted before that it generates a contraction semigroup. 
In fact, the semigroup $\mathcal{T}_1(t)$ generated by $A-M$ satisfies 
$\left|\left|\mathcal{T}_1(t)\right|\right|_{\mathcal{X}}=1$ for every $t>0$, hence $0=\omega_0(A-M)=s(A-M)$.
On the other hand we have $\left|\left|\exp\left\{t F'_{{\bf u}_*}\right\}\right|\right|\le\exp\left\{t\left|\left|F'_{{\bf u}_*}\right|\right|\right\}$, 
hence for the growth bound of the semigroup $\mathcal{T}_2(t)$ generated by $F'_{{\bf u}_*}$ we have $\omega_0\left(F'_{{\bf u}_*}\right)\le \left|\left|F'_{{\bf u}_*}\right|\right|$. 
We also note that $M+\nu\mathcal{I}$ is a dissipative operator and it generates a contraction semigroup, 
hence for the growth bound of the semigroup $\mathcal{T}_3(t)$ generated by $M$ we have $\omega_0(M)<-\nu$, 
and $\left|\left|\mathcal{T}_3(t)\right|\right|\le \exp\left\{-\nu\,t\right\}$. 
Finally, by applying a version of the Trotter product formula (see e.g. \cite[Corollary 5.8]{NAG}) we obtain 
$\omega_0(A-M+F'_{{\bf u}_*}+M)=\omega_0(A+F'_{{\bf u}_*})<0$ which shows that the non-trivial steady state ${\bf u}_*$ is locally asymptotically stable.   \hfill $\Box$

\begin{remark}
The stability condition in Theorem \ref{stability}  may seem very restrictive at the first glance. It may however be the case that 
the norm of the linearisation $\left|\left|F'_{{\bf u}_*}\right|\right|$ is small, especially since for our model 
we naturally have $\beta_0',b'_0,\beta'_2,b'_2\le 0$. On the other hand, even natural mortality (note that $\mu$ 
is natural mortality combined with infection induced mortality) is very high for example in the case of the mosquito 
{\it Aedes aegypti}, which is one of the target species 
for the introduction of {\it Wolbachia}, see \cite{McMeniman}.  We also note that the operator $-F'_{{\bf u}_*}$ 
may well be positive if the conditions of Theorem \ref{existencess} are satisfied, but there are no readily available 
results concerning the question that how does the infimum of the spectrum of an operator changes under positive perturbations. 
On the other hand if $F'_{{\bf u}_*}$ is dissipative and $\nu>0$ then similar arguments as used in the proof of 
Theorem \ref{stability} show that the steady state is locally asymptotically stable. 
\end{remark}

\begin{remark}\label{finalremark}
As we pointed out earlier, it is not possible to talk about linear sta\-bi\-li\-ty of the trivial steady state in the usual sense, as the nonlinearity $F$ is not 
Fr\'{e}chet differentiable at ${\bf 0}$. We established however global existence of solutions for any initial condition, and since condition (ii) in Theorem \ref{existencess} requires that the spectral bound of the linear part is positive in a small neighbourhood of ${\bf 0}$ 
it may be intuitively plausible to expect (for example utilising the variation formula \eqref{nonlinear2}) that the trivial steady state is unstable and the positive steady state is actually globally asymptotically stable 
if the condition of Theorem \ref{stability} holds true. 
\end{remark}

\section{Concluding remarks}

In this paper we introduced and analysed a nonlinear structured population model with diffusion in the state space.  
Individuals in the population are structured with respect to infection (for example bacterium) load, hence we used Wentzell boundary condition 
at the uninfected state $x=0$. Our model primarily intended to describe the evolution of {\it Wolbachia} infection in an arthropod, for example mosquito 
population. {\it Wolbachia} is a reproductive parasite and it affects the reproductive mechanisms of its host in an intriguing fashion. 
Here we focused on a cytoplasmic incompatibility (CI) inducing strain. Following \cite{Bourtzis,Guillemaud} we adopted the view of partial CI, 
namely, that a female can produce viable offspring only when mating with a male who has lower infection load. Therefore the functions 
$\beta_2,b_2$ are assumed to be monotone decreasing.  Note however, that these assumptions are not necessary to establish any of our results, 
such as Theorem \ref{existencesol} and Theorem \ref{existencess}. In fact, in our main result Theorem \ref{existencess}, the crucial assumption is 
the one which concerns the strict monotonocity of the functions $\beta_0,b_0$. The necessity of condition (i) in Theorem \ref{existencess} 
is also in agreement with the fact that CI itself does not regulate 
population growth, it only provides infected individuals with a reproductive advantage, see e.g. \cite{Engelstatter,FHin2,Sinkins}.  
The main difficulty in the mathematical analysis of our model also arises from this crucial assumption of the density dependent CI. 
 In particular, physiologically structured population models with infinite dimensional nonlinearities and with distributed recruitment terms, such as the one used here, are usually notoriously difficult to analyse, see e.g. \cite{CS2,FHin4}. 
The special form of the nonlinearity naturally appearing in $\beta_2$ and $b_2$ also implies that although the nonlinear operator $F$ is locally Lipschitz continuous at ${\bf 0}$ it is not Fr\'{e}chet differentiable. This also implies that even though intuitively one may expect that if the conditions of Theorem \ref{existencess} and Theorem \ref{stability} hold true then the trivial steady state is unstable, we cannot discuss linear stability of the trivial steady state using standard stability results from semilinear theory as developed for example in \cite{Henry,WEB}.

\section*{Acknowledgments}
\`{A}. Calsina was partially supported by the research projects DGI MTM2008-06349-C03-03 and 2009SGR-345.
J.~Z.~Farkas was supported by a University of Stirling Research and Enterprise Support Grant
and a Royal Society of Edinburgh International Travel Grant.

\end{document}